\input amstex
\documentstyle{amsppt}
\NoBlackBoxes
\loadmsbm
\loadeufm
\UseAMSsymbols

\NoRunningHeads

\topmatter

\title
The genus of curves on the three dimensional quadric
\endtitle

\author
Mark Andrea A. de Cataldo 
\endauthor

\address
Department of Mathematics, University of Notre Dame; Notre 
Dame 46556 (In), U.S.A.
\endaddress

\email
mark.a.decataldo.1\@nd.edu
\endemail



\abstract
By means of an {\it ad hoc} modification of the so-called 
``Castelnuovo-Harris
analysis" we derive an upper bound for the genus of integral
curves on the three dimensional nonsingular quadric which lie
on an integral surface of degree $2k$, as a function of 
$k$ and the degree $d$ of the curve. In order to obtain this 
we revisit the Uniform Position Principle to make its use 
computation-free. The curves which achieve this bound can be
conveniently characterized. 
\endabstract

\endtopmatter

\document

{\bf INTRODUCTION.}
\bigskip
The objects of investigation of this paper are the following 
two connected problems.
{\it What are the possible geometric genera of integral curves $C$ of 
degree $d$
lying on a nonsingular three dimensional quadric $Q_3$ in  $\Bbb P^4$ 
and on an integral surface $S$ of degree $2k$ contained in $Q_3$?}
As it is shown in this paper the above genera are bounded above by a 
function
of $d$ and $k$.
{\it What is the structure of the curves for which the genus is maximum
with respect to $k$ and $d$?}

The above problems are natural questions stemming from the analogue 
problems
that one can state by replacing, in what above, 
$Q_3$ by
$\Bbb P^3$ 
and $2k$ by $k$. These were answered completely in the paper
 [JH]. 
The paper [G-P] (and its refinement contained in  [E-P])
 deals with the very similar questions
 of {{\it (i)} determining
the biggest possible genus
for curves  of degree $d$ in 
$\Bbb P^3$
which do not lie on a surface of degree less than 
$k$, or lie on a surface of degree $k$
 and of {\it (ii)}  understanding the curves  for which 
the genus
is the maximum possible.
Going back to the quadric body ([A-S], \S6) gives an answer
to the problem of determining the maximum possible genus
for curves which lie on a surface of degree $2k$   under the assumption
$d>2k(k-1)$. To do so they use the technique of [G-P], coupled with 
the idea
of considering only hyperplane sections which are tangent to 
the quadric $Q_3$.

In this work an upper bound for the above genera is worked out
with no assumptions on the degree $d$. The bound is obtained
pursuing some numerical properties of embedded curves; a certain
maximization process is involved (cf. \S 2). 
 In analyzing the curves that should achieve that bound,
the unpleasant answer is that some systems of invariants are inconsistent
with each other so that, except for some special cases in which
the bound is sharp and the  curves of maximal possible genus
 are characterized, the bound turns
out to be {\it not} sharp: the biggest possible genus is
strictly smaller than the derived upper bound. 
It is
appropriate to say that some geometric information gets lost 
in the process.
At present the author is unable to bridge the gap
between the bound obtained in  this paper and ``the real bound."  
He 
conjectures that the extremal curves should be special curves 
(cf. Definition
1.3) so that the  right bound should then be (3.5.1) 
(see also \S 4, Question A).

The Paper [C-C-D] deals with questions  {\it(i)} and {\it (ii)} 
above in the context of
the arithmetic genus for curves in ${\Bbb P}^4$
of degree   ``sufficiently big."

The paper is organized as follows. \S 1 contains the
statement of Theorem 1.4, which is the main result of this paper:
it gives the upper bound, it says exactly when it is sharp
and characterizes the curves of maximum genus  in 
the cases in which the bound is sharp;
the proof of (1.4.1) is in \S 2 and  the one of (1.4.2) is in \S 3. 
This section also contains some preparatory material. The reader 
acquainted with the paper [JH] will realize how big is the debt of the 
present work towards it. However certain subtleties associated with the 
 possibility of having to deal with ``unbalanced"  curves
on nonsingular quadric surfaces had to be circumvented
by means of a systematic
use of the Uniform Position Principle; this section contains
some alternative formulations of it. \S 2 is a Castelnuovo-Harris type
approach to the determination of the wanted upper bound. 
It is an {\it ad hoc}
modification of the above mentioned paper of Harris. It contains also
a bound for the genus of curves for which the general hyperplane 
section is not contained
in any curve ``of type $k$" (see Theorem 2.10).
\S 3 discusses the bound obtained in \S 2. Moreover it deals with 
a special class
of curves (see Definition 1.3) which arise naturally 
in the contest of  curves with the biggest possible genus. 
 \S 4 is speculative in nature: it raises 
two questions 
that the author could not answer.

\bigskip
\noindent
{\bf Acknowledgments.} It is a pleasure  for the  author
to acknowledge crucial conversations with J. Migliore.
C. Peterson has kindly explained  folklore about Liaison.
Special thanks to the author's thesis advisor A.J. Sommese
for his guidance and his kind patience.
\noindent
This work was partially supported by a
``Borsa di studio per l'estero," n. 203.01.59 of the  C.N.R. of 
the Italian Government.
\bigskip
\vskip6truemm
{\bf 1. PRELIMINARIES.}\bigskip
The basic notation is the one of [Ha]. 

\noindent
The ground field is the field of complex numbers
$\Bbb C$. 

\noindent
$Q_i$ denotes a  smooth $i$-dimensional quadric in 
 a projective space $\Bbb P^{i+1}$.

\noindent
When there is no danger of confusion,
little  distinction is made between
Cartier divisors and associated rank one locally free sheaves and
the additive and tensor product notation are sometimes used
at the same time. The topological space will be sometimes dropped
when one is dealing with cohomology groups and their dimensions.

\noindent
In this paper the use of the adjective {\it general}  in connection
with an element $H$ of $\check\Bbb P$ 
is a quantifier; it means that there exists
a Zariski dense open subset $W$ of 
$\check\Bbb P,$ such that  for every
$H\in W$, ... 

\noindent
$\lfloor t\rfloor $
denotes the biggest integer smaller than or equal to 
$t$.

\bigskip

\bigskip
The following two sets of data are fixed throughout the  sequel of 
the paper:

\bigskip
\noindent
{\bf (1.1)}
{\it $C$ is an integral curve lying on a smooth three-dimensional quadric
$Q_3$, 
$k$ is  a positive integer,
$S_k$ is   an integral surface in
$|\Cal O_{Q_3}(k)|$
containing
$C$, 
$d$ 
and
$g$
 are the degree and the geometric genus of
$C$, respectively.}
\bigskip
\noindent
{\bf (1.2) Definition.} Define
$n_0$ and   $\epsilon$
when 
$d>2k(k-1)$
and
$\ \theta_0$ 
and 
$\epsilon'$
when
$d\leq 2k(k-1)$ 
as follows:

\noindent

$$\aligned 
&n_0:=\lfloor\frac{d-1}{2k}\rfloor+1;\\
d\equiv -\epsilon \quad &(mod \ 2k), \quad 0\leq \epsilon\leq 2k-1;\\
&\theta_0:=\lfloor \frac{d-1}{2k}\rfloor +1;\\
d\equiv-\epsilon'\quad &(mod\ 2\theta_0),\quad 0\leq
\epsilon'\leq 2\theta_0-1.
\endaligned
$$

\bigskip

The following class of curves plays a central role in the understanding 
of the curves whose genus is the maximum possible. Arithmetically
Cohen-Macauley is denoted by a.C.M..
\bigskip
\noindent
{\bf (1.3) Definition.} A curve $C$ as in (1.1) is said to be in the class
$\frak S (d,k)$,
if it is nonsingular, projectively normal and linked, in a complete 
intersection
on 
$Q_3$ 
of type 
$(k, n_0)$ 
if 
$d>2k(k-1)$  
($(\theta_0, k)$ 
if 
$d\leq 2k(k-1)$), 
to an  ({\it a fortiori}) a.C.M. curve 
$D_{\epsilon}$
(
$D_{\epsilon'}$)
of  degree
$\epsilon$
($\epsilon'$ 
respectively)
lying on a quadric surface 
hyperplane section of 
$Q_3.$

\bigskip
The following is the main result of this paper: it is a bound for 
the geometric genus of curves as in (1.1) in terms of
$d$ and $k$.
\bigskip
\noindent
{\bf (1.4) Theorem.}
{\it  Notation as in {\rm (1.1)} and {\rm (1.2)}. Assume 
$d>2k(k-1)$. Then

$$
g-1\leq \pi(d,k) -\Xi,\tag1.4.1
$$

where

$$
\pi(d,k)=
\left\{
\aligned
&\frac{d^2}{4k}+\frac{1}{2}(k-3)d-
\frac{\epsilon^2}{4k}-\epsilon(\frac{k-\epsilon}{2}),
\quad\qquad \qquad \qquad  0\leq \epsilon\leq k,\\
&\frac{d^2}{4k}+\frac{1}{2}(k-3)d-(k-\tilde \epsilon)
(\frac{\tilde \epsilon}{2}-\frac{\tilde \epsilon}{4k} + \frac{1}{4}),
\ \qquad  k+1\leq  \epsilon \leq 2k-1,\  \tilde \epsilon:=\epsilon-k; 
\endaligned
\right.
$$

and

$$
\Xi=\Xi(d,k)=
\left\{
\aligned
&0 \qquad \quad if \ \epsilon=0,\ 1,\ 2,\  2k-1,  \\
&1 \qquad \quad if \ else.
\endaligned
\right.
$$

\noindent
{\rm (1.4.2)} The bound is sharp for
$\epsilon=0$, $1,$ $2,$ $3,$ $2k-2,$ $2k-1$.
A curve achieves such a maximum possible genus if and only if
it is in the class
$\frak S (d,k)$, except, possibly, the cases 
$\epsilon=3$,  $2k-2$.

\noindent
Assume $d\leq 2k(k-1)$. Then the analogous statements with
$\pi'(d,k)=\pi(d, \lfloor \frac{d-1}{2k}\rfloor +1)= 
\pi(d, \theta_0)$ and with
 $\Xi',$
 $\epsilon'$  
$(\theta_0,k)$
and
$D_{\epsilon'}$
replacing
$\Xi,$ 
$\epsilon,$
$(k,n_0)$
and
$D_{\epsilon}$
 respectively, hold.}

\bigskip

The following, which is proven in [JH], page 194, is stated for 
the reader's
convenience; it is one of the two main ingredients of the analysis:
\bigskip 

\noindent
{\bf (1.5) Lemma} (Gieseker){\bf .} {\it Let 
$E\subseteq H^0({\Bbb P^1},{\Cal O}(l-1))$,
$\lbrace 0\rbrace \neq F\subseteq  H^0(\Bbb P^1,{\Cal O}(l))$
be two vector spaces of dimensions
$e$
and 
$f$
respectively, such that:
$E\times H^0(\Bbb P^1,{\Cal O}(1))$$\subseteq F$.
Then either 
$f\geq e+2$,
or
$|F|$
equals the complete linear system
$|{\Cal O}_{\Bbb P^1}(f-1)|$
plus
$(l-f+1)$ fixed points.}\bigskip

 The following two lemmata are nothing else but a reformulation 
 of the {\it Uniform Position Principle} (U.P.P.) (cf. [A-C-G-H], 
pages 111-113) in terms of subvarieties and 
of coherent sheaves respectively, rather than 
in terms of linear systems.
The use of this principle is the second main ingredient.
 First some  notation.\bigskip

\noindent

Let
 $\Cal C$
 be an integral curve 
of degree 
$d$
 in  a projective space 
${\Bbb P}$
 of any dimension, 
 $H$
 a hyperplane, 
$\Gamma$
the corresponding hyperplane section of 
$\Cal C$.

\noindent
 Let 
$\frak J$
be the incidence correspondence in 
$\Bbb P \times \check {\Bbb P}$
defined by
$\lbrace \,(p\ ;H)\mid  p\in H\,\rbrace$
 with first and second
  second projections  
$p$ and $q$ respectively, $\Cal F$ 
a coherent sheaf on 
$\Bbb P\times \check \Bbb P$. By abuse of notation  $H$  
can and will denote
the hyperplane and the corresponding point of
$\check\Bbb P$.

\noindent  
Let 
${\frak I}(\delta) $, 
$1 \leq \delta \leq d$
be the incidence correspondences in 
${\Cal C}^{\delta}\times \check {\Bbb P}$,
defined by
$\lbrace \,(p_1,\ldots, p_{\delta};H)\mid p_i\in H, \ \forall i\,\rbrace$,
where
${\Cal C}^{\delta}$
denotes the
$\delta$-fold product of 
${\Cal C}$.
The essence of the U.P.P. is that {\it the spaces
$\frak I(\delta)$
are irreducible}. This  principle should be regarded as a 
fundamental property
of curves in projective space.

\noindent
Finally define
$\hat\frak I(\delta)$
to be the quotient of 
$\frak I(\delta)$ 
by the action of the symmetric group
$S_\delta$: 
$\hat\frak I(\delta)=\frak I(\delta)/S_\delta$.
The spaces
$\hat\frak I(\delta)$
are irreducible as well.

\bigskip
\noindent
{\bf (1.6) Lemma} (U.P.P.1){\bf .} {\it Notation as above. Let
$\frak B$ 
be a closed subscheme of 
$\frak J$, 
${\frak B}_{H}\subseteq \frak B$
the closed subscheme cut by a general hyperplane
$H,$ i.e. 
$\frak B\cap q^{-1}(H)$. Then, either
$\Gamma\subseteq {\frak B}_H$,
or 
$\Gamma \cap {\frak B}_H=\emptyset$.}
\bigskip
\noindent
{\it Proof.} Define
$\Gamma':=\Gamma \cap {\frak B}_H$, and let 
$\delta$
be the cardinality of 
$\Gamma';$ $\delta$ is constant on a Zariski dense open subset
of
$\check{\Bbb P}$. 
Without loss of generality
assume
$\delta>0$.
  By shrinking the above set to another
Zariski dense open set
$W,$ if necessary, one can assume
that the incidence correspondence
$\frak I(\delta)$, restricted over
$W$, 
is a connected 
\'etale covering of degree
$\binom d\delta \delta!$.
Clearly the corresponding covering associated with 
$\hat\frak I(\delta)$ 
has degree
$\binom d\delta$.
 The assignment
$ \lbrace W \ni H\rbrace\longrightarrow \lbrace {\hat \Gamma'}
\in \hat \frak I (\delta)\rbrace$,
defines a holomorphic section over 
$W$ 
of the latter  covering; this is a contradiction unless
$\delta= d$. 
\qquad $\Cal Q.\Cal E.\Cal D.$

\bigskip
Let 
$B$
be any algebraic scheme of dimension
$b$. 
Consider the following decomposition:
$B_{red}= B_b\cup B_{b-1}\cup \ldots B_1\cup B_0$,
where
$B_i$ denotes the union of all the components of 
$B$ of dimension $i$ taken  with the reduced structure.

\bigskip
\noindent
{\bf (1.7) Lemma} (U.P.P.2){\bf .} {\it Notation as above. 
Let $H$ be a general 
hyperplane; consider the natural evaluation map:
$H^0({\Cal F}_{|H}) \otimes_{{\Cal O}_H} \check{\Cal F}_{|H}
\overset\eta\to\longrightarrow
{\Cal O}_H, $
 and 
${\Cal I}_{{\frak B}_H,H}:=Im (\eta)$. 
Let 
${\frak B'}_H$ be any closed subscheme supported at 
some 
$({\frak B}_H)_i.$
Then, either 
$\Gamma \subseteq\frak B'_H,$
or
$\Gamma \cap {{\frak B'}_H}=\emptyset.$} 

\bigskip
\noindent
{\it Proof.}
Generic flatness 
(cf. [Mu], Lecture 8) and 
semicontinuity   give
 a Zariski dense open subset
$W\subseteq \check {\Bbb P}$
over which
$q_*{\Cal F}$ 
is a locally free  coherent sheaf and
the natural maps
$q_*{\Cal F}\otimes_{{\Cal O}_W}k(w)\to h^0({\Cal F}_{|q^{-1}(w)}) $
are isomorphisms $\forall w\in W$.
Pick $\varsigma\gg 0$ such that 
$q_*\Cal F\otimes {\Cal O}_{\check{\Bbb P}}(\varsigma)$
is spanned by global sections on $\check \Bbb P$; then the following 
diagram commutes
and has surjective  vertical arrows, 
$\forall H\in W$:
$$
\CD
H^0(\Cal F \otimes q^*{\Cal O}_{\check \Bbb P}(\varsigma)) \otimes \check 
{\Cal F}
@>>>
{\Cal O}_{\frak J} \\
@VVV          @VVV  \\
H^0({\Cal F}_{|H}) \otimes_{{\Cal O}_H} \check {\Cal F}_{|H}  @>>>  
{\Cal O}_H.
\endCD
$$
\noindent
Hence the  scheme 
$({\frak B'}_H)_{red}$ 
is the restriction to $H$ of a closed subscheme
$\frak B'$  
in 
$\frak J.$
Let
$\Gamma'_{H}={\frak B'}_H\cap \Gamma$ 
and 
$\delta$ its cardinality; 
shrink $W$, if necessary, in order for 
$\delta$ to be constant over $W$. One can now conclude as 
in the previous lemma.
$\qquad \Cal Q. \Cal E. \Cal D.$

\bigskip
\noindent
 {\bf (1.8) Remark.} The above proposition is still valid, 
after obvious changes, 
if one replaces 
$\Bbb P$
 by some closed subscheme 
$ \Cal C\subseteq \frak T\subseteq \Bbb P$. 
In this paper
$\frak T=Q_3$.

\bigskip
\noindent
{\bf (1.9) Remark.} It is maybe worthy to observe that
(1.6) and (1.7) are both equivalent to the irreducibility
of the varieties
$\hat\frak I (\delta)$, $1\leq \delta\leq d$.

\bigskip

\vskip 6truemm
{\bf 2. DERIVING THE UPPER BOUND.}
\bigskip

The following is a  presentation of the relevant invariants
and of how to use them to give an upper bound on 
$g$
  as a function of
$d$ and $k$ (cf. [JH]).
\bigskip

Consider the following natural morphisms:
$\widehat C \overset\nu\to\rightarrow C 
\overset\iota\to\hookrightarrow Q_3 \hookrightarrow  {\Bbb P}^4$, 
where
$\widehat C \overset\nu\to\rightarrow C$
denotes the normalization of 
$C$,
 the other two arrows the given embeddings. 
All sheaves of the form
${\Cal O}(h)$ are pull-backs from
$\Bbb P^4$; the sheaves on 
$\hat C$
are pull-backs via 
$\nu$. Let
$\rho := \iota \circ \nu$ 
and
$\rho_l $
be the map induced in cohomology by
$\rho$;
 define:

$$\alpha_l:= \dim_{\Bbb C}[Im H^0( Q_3,{\Cal O}_{Q_3}(l)) 
\overset\rho_l\to\longrightarrow
 H^0(\widehat C, {\Cal O}_{\widehat C}(l)].
$$

\noindent
Let 
$H$ 
be a general hyperplane of
${\Bbb P}^4,$ 
$\Gamma:= C \cap H,$
$Q_2:= Q_3 \cap H$; 
then for every 
$l$ 
there is the map:
$H^0( Q_3, {\Cal I}_{\Gamma , Q_3} (l)) \overset\sigma_l\to\longrightarrow
H^0( \widehat C, {\Cal I}_{\Gamma , \widehat C}(l))\simeq H^0(\widehat C, 
{\Cal O}_{\widehat C} (l-1))$.

\noindent
Since
$Im (\rho_{l-1}) \subseteq Im(\sigma_l)$,
$Ker (\rho_l) = Ker (\sigma_l) = H^0({\Cal I}_{C,Q_3} (l))$,
$H^1_*({\Cal O}_{Q_3})=0$,
one gets the following chain of relations at the
end of which the quantities 
$\beta_l$ are defined:
$$
\aligned
\alpha_l - \alpha_{l-1}&=\dim\ Im(\rho_l) - \dim\ Im(\rho_{l-1}) \geq
\dim\ Im(\rho_l) - \dim\ Im(\sigma_l) \\
  &=
h^0({\Cal O}_{Q_3}(l)) - h^0({\Cal I}_{\Gamma , Q_3}(l)=
h^0({\Cal O}_{Q_2}(l))- h^0({\Cal I}_{\Gamma , Q_2}(l))\\
  &=:\beta_l.
\endaligned
$$
\noindent
One may think of 
$\beta_l$
as the number of independent conditions that
$\Gamma$
imposes on
$|{\Cal O}_{Q_2}(l)|$.

\noindent
Define:
$$
\gamma_l:=\beta_l - \beta_{l-1}=
[h^0({\Cal O}_{Q_2}(l)) - h^0({\Cal I}_{\Gamma , Q_2}(l))]-
[h^0({\Cal O}_{Q_2}(l-1)) - h^0({\Cal I}_{\Gamma , Q_2}(l-1))].
$$
\noindent
These ``second differences" are
quantities that can be  realized  geometrically as follows:
consider the following exact sequence defining a general conic 
$Q_1$ (to be chosen  so that it is smooth and it does not meet
$\Gamma$)
in 
$Q_2$:

$$ 0 \to {\Cal I}_{\Gamma, Q_2}(-1+l) \to {\Cal I}_{\Gamma ,Q_2}(l)
\to {\Cal  O}_{Q_1}(l) \to 0.$$

\noindent 
Let
$$
 E_l:= Im[ H^0({\Cal I}_{\Gamma, Q_2}(l)) \to H^0({\Cal O}_{Q_1}(l))],
$$
and
$$
e_l:= dim_{\Bbb C}(E_l).
$$ 
Then: 
$\gamma_l=h^0({\Cal O}_{Q_1}(l)) - [ h^0({\Cal I}_{\Gamma, Q_2}(l))-
h^0({\Cal I}_{\Gamma, Q_2}(l-1))]$. 

\noindent
It is now clear
that
$\gamma_l$
measures the incompleteness of the linear systems induced on
$Q_1$
by
$|{\Cal I}_{\Gamma , Q_2}(l)|$:
$$
\gamma_l= 2l+1-e_l.
$$
\noindent
Let:
$$
\theta:=\min\lbrace \,t\in 
{\Bbb N}\mid h^0({\Cal I}_{\Gamma,Q_2}(t)) >0\,\rbrace.
$$
\noindent
By the existence of
$S_k$,
one infers that
$\theta\leq k$.
\noindent
Let:
$$
n:=\min\lbrace \,\nu\in {\Bbb N}\mid \ |{\Cal I}_{\Gamma , Q_2}(\nu)| \ 
\   
\text{is  not empty
and  does not have  fixed  components}\, \rbrace .
$$
\noindent
Since 
$\gamma_l=\beta_l - \beta_{l-1},$
and 
$\beta_l=d,$
$\forall l\gg 0 $ 
(
$h^1({\Cal I}_{\Gamma, Q_2}(l)=0,$ 
$\forall l\gg 0 ),$
 one sees that
$\gamma_l=0,$ 
$\forall l\gg 0.$
 Define
$$
m:=\min \lbrace \,\mu \in {\Bbb N}\mid  \gamma_{\mu}=0 \,\rbrace;
$$
\noindent
Clearly
$\gamma_l=0$, $\forall l\geq m$.
\bigskip

 Following
Halphen, Castelnuovo, and more recently Gruson-Peskine and Harris,
by choosing 
$\lambda \gg 0$,
one gets:
$$
\align
g-1=  &\quad \hbox{\rm (Riemann-Roch)} \\
d\lambda -h^0({\Cal O}_{\widehat C}(\lambda)) {\leq } &
\quad (\alpha_l\leq h^0({\Cal O}_{\widehat C}(\lambda))) \\
d\lambda-\alpha_{\lambda}\leq    
&\quad ( \beta_t \leq \alpha_t - \alpha_{t-1})   \\
d{\lambda}- \sum_{t=0}^\lambda
\beta_t= &\quad (\beta_t=\sum_{l=0}^t\gamma_l) \tag2.1 \\
d\lambda-\sum_{l=0}^{\lambda}
(\lambda-l+1)\gamma_l= &\quad (\sum_{l=0}^{\lambda}\gamma_l=d) \\
  = &  \sum_{l=0}^\lambda(l-1)\gamma_l.
\endalign
$$

\noindent
The next step is to maximize the above sum with respect
to some  constraints on the numbers
$\gamma_l$. For the sake of clarity the analysis of these quantities
is divided into three cases: 
$d>2k(k-\ 1)$, 
$d\leq 2k(k-1)$ and
$\theta=\leq k-1$, 
$d\leq 2k(k-1)$ 
and
$\theta=k$. It is not necessary to distinguish between the
last two cases; however if one assumes
$\theta=k$
then one gets the smaller  upper bound (2.10), and does so without
assuming the existence of the surface
$S_k$.
\bigskip
{\bf The case:} 
$d>2k(k-1)$
{\bf.}
\bigskip

By   (1.2):
$d=2n_0 k-\epsilon$.

 If
$d>2k(k-1)$,
then equality holds in the inequality
$\theta\leq k$; 
for if 
one chooses 
$H$
 general then 
$D_k:=S_k\cap Q_2$ 
will be an integral curve which will not contain any of the components of
$D_{\theta}\in |{\Cal I}_{\Gamma, Q_2}(\theta)|$ 
so that, by computing intersections on 
$Q_2,$
 one gets: 
$D_k\cdot D_{\theta}=2\theta k\geq d>2k(k-1),$
that is
$\theta>(k-1)$.

\noindent
It follows that the linear systems
$E_l$ 
are empty in the range 
$[0,k -1]:$
$$
\gamma_l= 2l+1, \qquad \forall l \in [0,k -1].
$$

Since 
$D_k$
 is an integral curve
the linear systems
$|{\Cal I}_{\Gamma , Q_2} (l)|=D_k+|{\Cal O}_{Q_2}(l-k)|$, in the range
$[k,n-1]$,
so that, on the general 
$Q_1$, 
$E_l=D_k\cap Q_1 + |{\Cal O}_{Q_1}(l-k)|$;
it follows that:
$$
\gamma_i=2k, 
\qquad \forall i\in [k,n -1].
$$

The above interval is empty if and only if 
$|{\Cal I}_{\Gamma , Q_2} (k)|$
is free of fixed components, which in turn is equivalent 
to the statement that
$D_k$
 moves; this last condition implies of course
$h^0({\Cal I}_{\Gamma , Q_2} (k))\geq 2$
so that, if
$n=k$ 
then
$\gamma_k\leq 2k-1$.

As in [JH] it is now time to use Gieseker's
Lemma;
it allows to understand better the behavior of the quantities
$\gamma_l$ in the third remaining interval
$[n, m]$.
\bigskip
\noindent
{\bf (2.2) Lemma.} {\it
If
$k <n$ 
then one has the following information as to the behavior
of the quantities
$\gamma_l$: 
$\gamma_{n-1}-\gamma_{n}\geq 1$; 
$\gamma_{l-1} - \gamma_{l}\geq 2$, $\forall l \in [n, m-1]$;
$\gamma_{m-1}-\gamma_m=\gamma_{m-1}\geq 1$.
If
$k=n$ 
then  the same conditions hold except, possibly, the first one; in any case
$\gamma_k\leq 2k-1$}
\bigskip
\noindent
{\it Proof}. The only difference between the two
possibilities 
$k<n$
and 
$k=n$ lies, possibly,  in  
$(\gamma_{n -1}-\gamma_{n}).$
By what has been  shown above, the second statement for the case
$k=n$ is clear.

\noindent Assume therefore that 
$k<n$. 
One has
$E_j\subseteq H^0({\Cal O}_{Q_1}(j))\simeq H^0(
{\Cal O}_{{\Bbb P}^1}(2j))$
and
$H^0({\Cal O}_{{\Bbb P}^1}(1))\times 
H^0({\Cal O}_{{\Bbb P}^1}(1))\times E_{j-1}
\subseteq E_j$.
 One applies  Lemma (1.5) twice for every index 
$j$
 in the range considered, keeping in mind that,
since
$Q_1$
does not meet
$\Gamma$,
 the lack of fixed components
for 
$|{\Cal I}_{\Gamma , Q_2}(j)|$
implies the base-point-freeness of the corresponding
$E_j.$ It follows that
$e_l-e_{l-1}\geq 4 ,$ 
except possibly 
$l=n,$ $m$
where
$e_l-e_{l-1}\geq 3.$ \qquad
$\Cal Q.\Cal E. \Cal D.$

\bigskip
Since 
$|{\Cal I}_{\Gamma , Q_2}(n)|$
does not have fixed components, any curve in that linear system
cuts on 
$D_k$
 a set of 
$2n k$
 points (counted  with multiplicities) that contains
$\Gamma,$
so that
 $2n k\geq d:$
$$
n\geq n_0=\lfloor \frac{d-1}{2k}\rfloor+1.
$$

One can summarize the information on 
$\gamma : [0,m]\to {\Bbb N}$ 
as follows:
$$
\alignat 2
\gamma_l &=2l+1, &\qquad &l\in [0,k-1];  \\
\gamma_l &=2k, &\qquad & l\in [k,n-1];  \\
\gamma_{n} &\leq 2k-1;&&   \\
\gamma_{l}-\gamma_{l+1} &\geq 2,&\qquad & l\in [n, m-2];   \\ 
\gamma_l&=0,&\qquad &l\geq m;\\
\sum_{l=0}^m \gamma_l&=d. & & 
\endalignat
$$

After (2.1), the goal is to maximize 
$\sum_{l=0}^m(l-1)\gamma_l$,
subject to the above constraints.
One can start by reducing  the process to the case in which
$n=n_0$. 
\bigskip
\noindent
{\bf Remark.} It should be noted that 
$m\leq n_0+k$.
 This is a
 straightforward consequence of the constraints on   
$\gamma$. 
In particular
one could already find an a priori upper bound for 
$\sum (l-1)\gamma_l$
 by
adding up setting, for example, 
$\gamma_l=2k$.
\bigskip

\bigskip
\noindent
{\bf (2.3) Lemma.} {\it Given any function
$\gamma$ 
subject to the above constraints there exists
a function
$\tilde\gamma$,
subject to the same constraints, for which the corresponding
$n=n_0$ 
(here 
$n$
 is the first number greater or equal to
$k$
 for which 
$\gamma_{n}<2k)$
 and for which 
$\sum (l-1)\gamma_l\leq \sum (l-1)\tilde \gamma_l$.} 

\bigskip
\noindent
{\it Proof.} Assume 
$n -n_0=:\xi>0$,
otherwise there is nothing to show.  One has:
$$
d=\sum_{0}^{n-1}\gamma_l +\sum_{n}^m\gamma_l=
 -k^2 + 2n_0k + 2\xi k + \sum_{n}^m\gamma_l=
-k^2 + d + \epsilon +2\xi k +\sum_{n}^m\gamma_l;
$$
\noindent
it follows that
$$
k^2=\epsilon + 2\xi k + \sum_{n}^m\gamma_l.
$$

\noindent
By the above
$k\geq 2$,
so that 
$2\xi k\geq 4$
 and  
$\sum_{n}^m\gamma_i\leq k^2-4$.
It follows that one of the 
following conditions must hold:

\noindent
{\it a)} there is an index 
$n\leq j \leq m-1$
 for which 
$\gamma_{j -1}-\gamma_{j}\geq 5$;

\noindent
{\it b)} there are  two distinct indices
 $j_1<j_2$
 as in {\it a)} for which
$\sum_1^2(\gamma_{j_t-1}-\gamma_{j_t})\geq 6$;

\noindent
{\it c)} there are three distinct indices 
$j_1<j_2<j_3$,
as in {\it a)}
such that 
$\sum_1^3(\gamma_{j_t -1}-\gamma_{j_t})\geq 8$;

\noindent
{\it d)} there are four distinct indices
$j_1<j_2<j_3<j_4$, 
as in {\it a)} such that 
$\sum_1^4(\gamma_{j_t -1}-\gamma_{j_t})\geq 10$.

\noindent
In case {\it a)} one decreases (increases) 
$\gamma_{j-1}$ 
($\gamma_j$)
by one. In case {\it b)} either one is also
in case {\it a)} or one can decrease (increase) 
$\gamma_{j_1-1}$ ($\gamma_{j_2}$) by one. Similarly in the remaining
cases.
As a consequence of this process,  the constraints are respected but 
$\sum_0^m(l-1)\gamma_l$
increases. Since this sum is bounded from above by the above remark
the process must come to an end, {\it i.e.} one can modify any
$\gamma$
 to a 
$\tilde \gamma$ 
for which the corresponding
$\xi=0.$ 
$\qquad \Cal Q.\Cal E.\Cal D.$

\bigskip
\noindent
{\bf (2.4) Corollary.} {\it The following function 
$\tilde \gamma$
satisfies the constraints
and maximizes 
$\sum_0^m(l-1)\gamma_l:$}

{\it if
$0\leq \epsilon\leq k$:}
$$
\aligned
\tilde \gamma_l =& 2l+1, \\
\tilde \gamma_l =& 2k,    \\
\tilde \gamma_l =&2(k+n_0-l)-1,  \\
\tilde \gamma_l = &[2(k+n_0-l)-1]-1,\qquad  \qquad\\
\tilde \gamma_l =&0, \endaligned
\aligned 0&\leq l \leq k-1, \\
k&\leq l<n_0, \\
n_0&\leq l\leq n_0+k-\epsilon -1, \\
n_0+k-\epsilon&\leq l\leq n_0+k-1, \\
n_0+k&\leq l,
\endaligned
$$

{\it if 
$k+1\leq \epsilon\leq 2k-1$, 
let
$\tau:=\epsilon -k$,
then:}

$$
\aligned
\tilde \gamma_l =& 2l+1, \\
\tilde \gamma_l =& 2k,    \\
\tilde \gamma_l =&[2(k+n_0-l)-1]-1,  \\
\tilde \gamma_l = &[2(k+n_0-l)-1]-2,\qquad  \qquad\\
\tilde \gamma_l =&0, \endaligned
\aligned 0&\leq l \leq k-1, \\
k&\leq l<n_0, \\
n_0&\leq l\leq n_0+k-\tau -1, \\
n_0+k-\tau&\leq l\leq n_0+k-1, \\
n_0+k&\leq l,
\endaligned
$$

\bigskip
\noindent
{\it Proof.} By the previous lemma one can assume 
$n =n_0$;
it remains to define  
$\tilde \gamma_l$
in such a way  that
$\sum_{n_0}^{m} \gamma_l = k^2-\epsilon$
 and 
$\sum(l-1)\gamma_l $
is maximized. First define
$\breve\gamma_l=2(k+n_0-l)-1$
for 
$l\in [n_0, n_0 +k -1]$.
Now one has to delete from the graph of 
$\breve\gamma$
$\epsilon$
 points; 
$\tilde\gamma$
 is the way to delete those points while
 maintaining  the constraints and meeting the 
 above maximization requirements.
$\qquad \Cal Q. \Cal E. \Cal D.$
\bigskip

\noindent
{\bf (2.5)}\quad If one adds up 
$\sum_{l=0}^{n_0+k}(l-1)\tilde{\gamma_l}$,
one gets the desired function
$\pi=\pi(d,k)$ for which
$(g-1)\leq \pi$. For its explicit form see Theorem 1.4.

\bigskip
\noindent
{\bf  Remark.} The above is the bound obtained in [A-S], \S6 for curves
$C$
of degree
$d>2k(k-1)$ 
 contained in an integral surface of degree 
$2k$.
As it will be shown in \S 3, the bound (2.5)
is not quite sharp.

\bigskip

\vskip 6truemm
{\bf  The case:} 
$d\leq 2k(k-1)$, $\theta\leq k-1${\bf .}
\bigskip

In this case the analysis of the behavior of the function 
$\gamma$
 associated with
$C$
 is   analogous to the first case. The twist
is the behavior of 
$\gamma$ 
in the interval
$[\theta , n -1]$. 
The following takes care of that interval.
\bigskip
\noindent
{\bf (2.6) Proposition.}  
$\gamma_l=2\theta$, 
$\forall l\in [\theta, n -1]$.
\bigskip
\noindent
{\it Proof.} Let 
$l$
 be in the above range. Using the notation of (1.7) define
$\frak J_{Q_3}:=p^{-1}
Q_3,$ 
$\frak C:=p^{-1}C$
 and define 
$\Cal F(l):={\Cal I}_{ \frak C, 
\frak J_{Q_3}}\otimes p^*{\Cal O}_{Q_3}(l).$ 
The proof of Lemma 1.7 and Remark 1.8 imply that 
for every $l$ 
the fixed component  
$F_l$ of 
$|{\Cal I}_{\Gamma , Q_2}(l)|$
contains all of
$\Gamma$.
Clearly
$F_\theta\supseteq F_{\theta +1}\supseteq\ldots \supseteq F_{n -1}$.

\noindent
If 
$F_\theta\supsetneq F_l$,
for some $l$, then the curve
$F_\theta - F_l$
would be free to move in
$|{\Cal I}_{\Gamma , Q_2}(\theta)|$, a contradiction. It follows that
$F_\theta=\ldots=F_{\varphi -1}$. 

\noindent
To conclude one has to show that 
$F_\theta$
is actually a member of
$|{\Cal I}_{\Gamma , Q_2}(\theta)|$.

\noindent
One can  choose a line
$\ell \subseteq \check{\Bbb P}^4$
such that:

{\it i)} it defines a pencil of hyperplane sections of
$Q_3$ based on a smooth conic 
$\bar Q_1$ 
that does not meet $C $,

{\it ii)} it meets the open set $W$ of (1.7)

and

{\it iii)} it meets the open set of $\check {\Bbb P}^4$ for which
$\Gamma$ has cardinality $d.$

\noindent
Using the same method  as in  the quoted lemma one constructs a surface
$\tilde S$ on 
$q^{-1}(\ell)$ (which is the blowing up of
$Q_3$ along 
$\bar Q_1$),
that cuts the general element of the pencil the corresponding curve
$F_\theta$. 
This surface descends to
$Q_3$ as a surface
$S$ 
that cuts on, the general  element of the pencil, a  curve of the form
$F_{\theta}+ \mu \bar{Q}_1$, 
where 
$\mu$ 
is some integer.
Since 
$Pic(Q_3)\simeq Pic(\Bbb P^4)$,
one sees that
$S\in |\Cal O_{Q_3}(\zeta)|$, for some integer 
$\zeta$;
 it follows that
$F_\theta \in |\Cal I_{\Gamma, Q_2}(\chi)|$, for some integer
$\chi$. By the minimality of $\theta$ one concludes 
$\theta = \chi$. 
$\qquad \Cal Q. \Cal E. \Cal D.$
\bigskip

\noindent
{\bf  Remark.}
The same method as above offers an alternative way to prove, 
less elementarily
but in an unifying way, that
$\gamma_l=2k,$ $\forall l\in [k, n -1]$ in the case
$d>2k(k-1).$

\bigskip
Now one can repeat the analysis of the case
$d>2k(k-1)$ and obtain an analogous function
$\tilde \gamma$
as follows:  substitute
$k$ and $n_0$ 
by
$\theta_0$ and $ k$ 
respectively,
 for if one does so
then 
$m$ 
will be maximized.
\bigskip

\noindent
{\bf (2.7)} Adding up one gets, as in (2.5), a function
$\pi'=\pi'(d,k)$
that bounds $g-1$
from above. By construction
$\pi'(d,k)=\pi(d, \lfloor \frac{d-1}{2k}\rfloor +1)=\pi(d,\theta_0)$.

\bigskip
 \noindent
 {\bf  Remark.} The bound $\pi'$ is not quite sharp as well 
(see \S 3.).

\bigskip

\vskip 6truemm
{\bf The case:} $d\leq 2k(k-1)$, $\theta=k${\bf .}

\bigskip
In what follows the surface $S_k$ will play no role. Hence the only 
assumptions needed are:

\bigskip \noindent
(2.8)\quad
{\it $C\subseteq Q_3$ 
is an integral curve of degree
$d\leq 2k(k-1)$,
for which the general  hyperplane section
$\Gamma \subseteq Q_2$  is not contained in any curve
belonging to the linear system
$|{\Cal O}_{Q_2}(k-1)|.$}

\bigskip
Clearly
$\theta\geq k$; as in the previous case
$\gamma_l=2\theta$, if 
$l\in[\theta, n-1]$; also Lemma 2.2 holds
with $k$ 
replaced by
$\theta$.

Now one starts modifying
$\gamma$,
 if necessary, to maximize
$\sum (l-1)\gamma_l$.
\noindent
First of all, since
$|{\Cal O}_{Q_2}(k-1)|\simeq {\Bbb P}^{k^2}$,
one has
$d>k^2=\sum_{l=0}^{k-1}\gamma_l$. Next, since the numbers
$\gamma_l$
must add up to
$d\leq 2k(k-1)$, 
after reducing oneself, as in Lemma 2.3,
to the case 
$k=\theta=n$, it is easy to see which function 
$\tilde \gamma$ 
 maximizes
$m$, and thus 
$\sum(l-1)\gamma_l$:

let 
$\nu, \epsilon$ 
 be the unique non-negative integers such that
$$
d=k^2+ {\nu}^2+\epsilon,\qquad 0\leq \epsilon\leq 2\nu,\tag2.9
$$

then define
$\tilde \gamma$
as follows:

if 
$0\leq \epsilon\leq \nu$, then
$$
\aligned
\tilde \gamma_l =& 2l+1, \\
\tilde \gamma_l =&[2(k+\nu-l)-1]+1,  \\
\tilde \gamma_l = &[2(k+\nu-l)-1],\qquad  \qquad\\
\tilde \gamma_l =&0, \endaligned
\aligned 0&\leq l \leq k-1; \\
k&\leq l\leq k+\epsilon -1; \\
k+\epsilon&\leq l\leq k+\nu -1; \\
k+\nu&\leq l;
\endaligned
$$
 if 
$\nu+1\leq \epsilon\leq 2\nu$,
let first
$\tau:=\epsilon-\nu$ and
$$
\aligned
\tilde \gamma_l =& 2l+1, \\
\tilde \gamma_l =&[2(k+\nu-l)-1]+2,  \\
\tilde \gamma_l = &[2(k+\nu-l)-1]+1,\qquad  \qquad\\
\tilde \gamma_l =&0, \endaligned
\aligned 0&\leq l \leq k-1; \\
k&\leq l\leq k+\tau -1; \\
k+\tau&\leq l\leq k+\nu -1; \\
k+\nu&\leq l.
\endaligned
$$
\bigskip 

\noindent
 {\bf  Remark.} Even without adding up, at this point
 one already knows, since
$\theta_0<k$,
 that the result
will be strictly smaller than the corresponding
$\pi'$ of (2.7).

\bigskip
The proof of (1.4.1) is now complete. By adding up what above 
one gets the following:

\bigskip \noindent
 {\bf (2.10) Theorem.} {\it Assumptions and notation as in {\rm (2.8)}
and {\rm (2.9)}. The geometric genus of $C$ satisfies
the following bound:}
$$
g-1\leq
\left\{
\aligned
&(k-\frac{3}{2})d-\frac{1}{3}(k^3-\nu^3)-\frac{1}{6}(k-\nu)
+\frac{1}{2}\epsilon^2,
\quad \quad \ \ \hbox{\rm if \ $0\leq \epsilon\leq \nu$;} \\
&(k-\frac{3}{2})d-\frac{1}{3}(k^3-\nu^3)-\frac{1}{6}(k-\nu)
+\frac{1}{2}\nu^2+\\
&+\frac{1}{2}(\epsilon -\nu)(\epsilon -\nu +k -3),
 \qquad\quad\qquad\  \ \ \quad \qquad \hbox{\rm if \ $\nu
+1 \leq \epsilon\leq 2\nu$.}
\endaligned
\right.
$$
\bigskip
\vskip 6truemm
{\bf 3. DISCUSSION: When is the bound sharp? When  is it not?}

\bigskip

Assume the curve 
$C$
has geometric genus maximum with respect to the upper bounds 
$\pi$, $\pi'$ of (2.5) and (2.7).
\noindent
 In particular
$\gamma=\tilde \gamma$ and the inequalities in (2.1) are all equalities.
By the following elementary claim, {\it if such a curve exists} 
then it will be
smooth and projectively normal.

\bigskip
\noindent
{\bf (3.1) Claim.}   
$C$
{\it is smooth if and only if} 
$\rho_l$
{\it is surjective}
$\forall l \gg 0$. 
{\it Moreover if} 
$C$ 
{\it is smooth it is projectively normal} (i.e.
$\rho_l$
is surjective
$\forall l$)
{\it if and only if} 
$\beta_l = \alpha_l - \alpha_{l-1},$ $\forall l.$
 
\bigskip
\noindent
{\it Proof.} (Cf. [JH], page 193). The first part is clear since
the normalization map 
$\nu^*:\Cal O_C\to \Cal O_{\hat C}$ has zero cokernel if and only if
$C$ is  smooth. As to the second part one argues as follows. If
$\rho_l$ is surjective for every $l$,
then
$\sigma_l=\rho_{l-1}$ for every $l$ as well.
Conversely  assume $C$
is not projectively normal and let 
$l_0$ be any  index such that
$\rho_{l_0+1 }$ 
is surjective but 
$\rho_{l_0}$ 
is not. 
Since $h^1(\Cal I_{C, Q_3}(l_
0+1))=0$, 
$\sigma_{l_0+1}$ 
is surjective. It follows that 
$\sigma_{l_0+1}> \rho_{l_0}$, so that
$\alpha_{l_0+1}-\alpha_{l_0}>\beta_{l_0}$. 
$\qquad \Cal Q.\Cal E. \Cal D.$

\bigskip
\noindent
{\bf (3.2) Claim.} {\it If 
$d>2k(k-1)$
there exists an integral surface
$S_{n_0}\in |{\Cal I}_{C,Q_3}(n_0)|$
such that  
$S_k\nsubseteq S_{n_0}$. If
$d\leq 2k(k-1)$ 
then there exists an integral surface 
$S_{\theta_0} \in |\Cal I_{C, Q_3}(\theta_0)|$.

\bigskip
\noindent Proof.} 
Assume first that 
$d>2d(k-1)$.
Then since
$|{\Cal I}_{\Gamma,Q_2}(n_0)|$
is free of fixed components, one finds  in it an element 
$ F_{n_0}$ 
that does not contain the irreducible curve 
$D_k$. The projective normality of 
$C$ 
translates into the surjection 
$H^0({\Cal O}_{Q_3}(l))\twoheadrightarrow H^0({\Cal O}_{C}(l))$,
$\forall l$. This, in turn, is equivalent to
$H^1({\Cal I}_{C,Q_3}(l))=0$,
$\forall l$. Applying this to the case 
$l=n_0-1$ one gets the surjection
$H^0({\Cal I}_{\Gamma, Q_3}(n_0))\twoheadrightarrow 
H^0({\Cal I}_{\Gamma, Q_2}(n_0))$.

\noindent
Therefore
$F_{n_0}$  
can be  lifted to a surface  
$S_{n_0}\in |\Cal I_{C,Q_3}(n_0)|$. 
Since 
$Pic( Q_3)=\Bbb Z$ 
it follows that this surface 
is integral otherwise one would find 
$n_1<n_0$
for which there is an element
$F_{n_1}\in |{\Cal I}_{\Gamma,Q_2}(n_1)|$
not containing
$D_k$,
 a contradiction, since then
$|{\Cal I}_{\Gamma,Q_2}(n_1)|$
would be free of fixed components.

\noindent
If 
$d\leq 2k(k-1)$ 
then there is a unique element 
$F_{\theta_0}\in  |{\Cal I}_{\Gamma, Q_2}(\theta_0)|$; one can lift it
to a surface 
$S_{\theta_0} \in |\Cal I_{C,Q_3}(\theta_0)|$ 
which
is integral by the minimality property of
$\theta_0$. \qquad $\Cal Q.\Cal E. \Cal D.$

\bigskip
\noindent
By what has just been shown, 
$C$ 
is residual to a curve 
$D_\epsilon$
of degree 
$\epsilon$ 
if 
$d>2k(k-1)$ 
($D_{\epsilon'}$ 
of degree
$\epsilon'$ 
if 
$d\leq 2k(k-1)$)
 in a complete intersection on 
$Q_3$ 
of type
$(k,n_0)$ 
($(\theta_0, k$), respectively).

\bigskip
\noindent 
The following lemma is the technical device
needed to relate
$C$ and 
$D_\epsilon$ 
($D_\epsilon'$).
 The proof is a mere
generalization of [JH], page 199, where the case
$S\simeq  \Bbb P^2$ was dealt with. It will be used here only in the case
$S\simeq Q_2$; proving it in a more general form 
is not more costly.

\bigskip
\noindent
 {\bf (3.3) Lemma.} {\it  Let
$S$
be a normal and projective surface, 
${\Cal O}_{S}(1)$ 
a nef and big line bundle on it, 
$F$ 
and
$G$ 
two curves in
$|{\Cal O}_{S}(n)|$
and
$|{\Cal O}_{S}(m)|$
respectively without any common component. Denote by
$\tilde \Gamma$
their scheme-theoretic intersection. 
Assume $\tilde \Gamma=\Gamma + \Gamma'$,
where
$\Gamma $
is reduced and disjoint from
$\Gamma'$
and
$\Gamma \subseteq  S_{reg}$. Then:
$$
h^1(S, {\Cal I}_{\Gamma ,S}(n+ m -l)\otimes\omega_{S})=
h^0(S,{\Cal I}_{\Gamma',S}(l)),
\qquad \forall l < m,\ n.
$$}

\noindent
{\it Proof.} Let
$\pi:S'\to S$
be the blowing up of 
$S$ 
along
$\Gamma$,
$E$
the exceptional divisor. Since 
$F$
and
$G$
meet transversally at 
$\Gamma$ 
one gets the following relations concerning  strict transforms:  
$F'=\pi^*F-E$,
$G'=\pi^*G -E$. 
Denote by 
$\Gamma''$
the scheme on 
$S'$
isomorphic to 
$\Gamma'$
via 
$\pi$, and
by
$\Cal O_{S'}(\upsilon)$ the pull back 
$\pi^*\Cal O_{S'}(\upsilon).$
 By taking the cohomology of the following resolution:
$$
0\to \pi^*{\Cal O}_{S}(-n-m+l)+2E\to \pi^*{\Cal O}_S(-n+l)+
E\oplus\pi^*{\Cal O}_S(-m+l)
+E\to {\Cal I}_{\Gamma''}(l)\to 0,
$$
one gets, for
$l<n,$ $m:$
$$
\align
0\to H^0( {\Cal I}_{\Gamma''}(l))
&\overset b\to\rightarrow H^1({\Cal O}_{S'}(-n-m+l)+2E)\to \\
 &H^1({\Cal O}_{S'}(-n+l)+E)\oplus H^1({\Cal O}_{S'}(-m+l)+E)\to\ldots
\endalign
$$
\noindent
The above vector space  is zero,
for $l<n,$ $m,$ as it is now shown.  Leray spectral sequence gives
$H^1({\Cal O}_{S'}(-t+l) +E)=H^1({\Cal O}_S(-t+l))$,
$\forall t$.
The latter group is zero
(this is a well-known argument): take a desingularization 
$\Cal S \to S$, 
 pull back 
${\Cal O}_{S}(-t+l)$ 
to a  nef and big
${\Cal O}_{\Cal S}(-t+l)$;
Kawamata-Viewheg vanishing (cf. C-K-M, Lecture 8) descends, again by 
Leray 
spectral sequence, to 
$S$.

\noindent 
Next, 
$S'$ being normal it is Cohen-Macauley. Using Serre Duality:
$$
H^1({\Cal O}_{S'}(-n-m+l)\otimes 
{\Cal O}_{S'}(2E))\simeq H^1(\pi^*(\omega_S 
\otimes{\Cal O}_{S}(n+m-l))\otimes{\Cal O}_{S'}
(-E))^{\vee}.
$$
By  Leray spectral sequence one concludes using the isomorphism 
$b$. 
$\qquad \Cal Q.\Cal E.\Cal D.$

\bigskip

\noindent
{\bf (3.4) Claim.} {\it 
$D_{\epsilon}$ 
lies on some quadric surface
$\Sigma\subseteq Q_3$.
\bigskip
\noindent
Proof.} If 
$\epsilon =0$ there is nothing to prove. Let 
$\epsilon>0$,
and denote by
$\Gamma'$
  the general hyperplane section of
$D_{\epsilon}$.
 Then 
$\gamma_{n_0+k-1}=0$ and 
$\beta_{n_0+k-2+t}=d, \forall t\geq 0$. 
One has:
$0<\gamma_{n_0+k -2}=d-\beta_{n_0+k-3}=
d-[d-h^1({\Cal I}_{\Gamma, Q_2}(n_0+k-3))]=h^0({\Cal I}_{\Gamma'}(1))$;
the last equality follows from (3.3). 
$C$
being projectively normal, 
$D_\epsilon$
is  a.C.M. (cf., for example, [Mi], Th 1.1).
It follows that 
$h^0({\Cal I}_{D_\epsilon,\  \Bbb P^4}(1)>0$, i.e. 
$D_{\epsilon}$
is contained in a hyperplane.
\qquad$\Cal Q. \Cal E. \Cal D.$
\bigskip
\noindent

\bigskip
\noindent
{\bf (3.5) Example.} Here the function 
$\gamma$ and the genus of the curves in the class $\frak S(d,k)$ 
are computed.
The curves in these classes are the natural candidates to be the curves
of maxima genera with respect to $d$ and $k$.

\noindent
Since the two cases 
$d>2k(k-1)$
and
$d\leq 2k(k-1)$
are treated in the same way, the example is worked out only 
in the former case. 

\noindent
Let 
$\Gamma'$ denote
the general hyperplane section of
$D_{\epsilon}$. Assume  first that  
$\epsilon$ 
is even: $\epsilon=2\alpha$. 
Since 
$D_{\epsilon}$ 
is a.C.M., one sees that 
$D_{\epsilon}\in |\Cal O_{\Sigma}(\alpha)|$.
One takes the
cohomology of the following  projective resolutions of the 
twists of the ideal 
sheaf of 
$\Gamma'$:
$$
0\to {\Cal O}_{Q_2}(-1-\alpha+l)\to {\Cal O}_{Q_2}(-1+l)\oplus
{\Cal O}_{Q_2}(-\alpha +l)\to{\Cal I}_{\Gamma' , Q_2}(l)\to 0.
$$
To compute the quantities 
$\gamma_l$
one argues as in (3.4) using Lemma 3.3: 
$\gamma_{n_0+k-2-l}=h^0({\Cal I}_{\Gamma',Q_2}(l+1))-
h^0({\Cal I}_{\Gamma',Q_2}(l))$, 
$\forall l\leq k-2$.

\noindent
Now it is assumed that 
$\epsilon $
is  odd:
$\epsilon=2\alpha -1$. One can pick a line
$\Cal L$ on 
$\Sigma$ so that
$\Cal M:=D_{\epsilon}\cup \Cal L$ is a curve in
$|\Cal O_{\Sigma}(\alpha)|$ (cf. [A-C-G-H], Ex. III D7). 
The general  hyperplane section of 
${\Cal M}$
is
$\Gamma''=\Gamma'\cup p$, where 
$p $ 
is the point hyperplane section of the line $\Cal L$.
In addition to the projective resolution for 
$\Gamma''$, which is the same as above, 
one also has the following exact sequences:
$$
0\to {\Cal I}_{\Gamma'',Q_2}(l) \to{\Cal I}_{\Gamma', Q_2}(l)\to 
{\Cal O}_p\to 0.
$$ 
\noindent
Keeping in mind that
$h^1({\Cal I}_{\Gamma'',Q_2}(l))=0,$ $\forall l\geq \alpha$, 
Lemma 2.2 and the usual constraints
 a straightforward computation,  analogous to the one of the case
$\epsilon $ 
even, gives the desired quantities 
$\gamma$.

\noindent
From what above  one concludes that the function 
$\hat\gamma$ 
for these
special curves is the following:

first let
$$
\Delta:=\left\{
\aligned
&0 \qquad\hbox{\rm if $\epsilon=0$ or $\epsilon$ is odd},\\
&1 \qquad\hbox{\rm if $\epsilon$ is even and $\epsilon\geq 2$};
\endaligned
\right.
$$
and let
 $\alpha$ be as above, then
$$
\alignat 2
\hat \gamma_l&= [2(n_0+k-l)-1],&\qquad & n_0\leq l\leq n_0+k-\alpha -2\\
\hat \gamma_l&=[2(n_0+k-l)-1)]-\Delta,
&\qquad & l=n_0+k-\alpha-1,\\
\hat\gamma_l&=[2(n_0+k-l)-1]-2, &\qquad & n_0+k-\alpha\leq l\leq n_0+k-2\\
\hat\gamma_l&=0,&\qquad &n+k-1\leq l.
\endalignat
$$
Morevover by adding up one gets that the genera of these curves are:
\bigskip
\noindent
(3.5.1)\centerline{$
g-1=\Pi:=\left\{
 \aligned
&\frac{1}{4k}d^2+ \frac{1}{2} (k-3)d -
\frac{\epsilon}{2}[(k-1)(1-\frac{\epsilon}{2k})]-\frac{1}{4},
\quad\ \   \hbox{\rm if $\epsilon$ is odd,} \\
&\frac{1}{4k}d^2+\frac{1}{2}(k-3)d -
\frac{\epsilon}{2} [(k-1)(1-\frac{\epsilon}{2k})],
\qquad\qquad \hbox{\rm if $\epsilon$ is even}.
\endaligned 
\right.$}

\bigskip
\noindent
The two functions 
$\tilde\gamma$ 
of (2.4) 
and 
$\hat \gamma$ 
coincide if and only if
$\epsilon=0,$ 
$1,$ 
$2,$ 
$2k-1$.
This proves that  the geometric genus of
$C$
achieves the bound 
$\pi$
 if and only if
$\epsilon=0,$
$1,$
$2,$
$2k-1$
and
$C\in \frak S (d,k)$.
 Conversely if $C\in \frak S (d,k)$ 
then its associated function
$\gamma$  and its genus are as in (3.5.1).

\bigskip \noindent
As to the cases
$\epsilon=3$, $2k-2$. By what above
$g-1<\pi$, and
$\sum (l-1)\hat {\gamma_l}=\pi-1$; it follows
that this latter value is the sharp bound.

\bigskip
This proves  (1.4.2) so that the proof of Theorem 1.4 is now complete.
\bigskip
\vskip 6truemm
{\bf 4. TWO OPEN QUESTIONS.}
\bigskip

The author would like to pose the  following two questions. The first one
is the  consequence of the
incompleteness of  Theorem 1.4.
The answer to the second one  would constitute
a natural property of curves on quadrics.

\bigskip
{\bf Question A.} {\it Is it true that the curves of maximal genus
with respect to $(d,k)$
are the ones of the class
$\frak S(d,k)$?}
\bigskip

A positive answer would give the sharp bound (3.5.1)
and the  complete characterization
of the curves of maximal genus.

\bigskip{\bf Question B.} {\it The  {\rm U.P.P.} is expressed in terms of
the space of hyperplane sections of $Q_3$, {\it i.e.}
$\check{\Bbb P}^4$. Does the analogue statement hold
if one considers only  hyperplanes which are tangent
to $Q_3$, {\it  i.e.} replacing
$\check{\Bbb P}^4$ by 
$\check{Q}_3$? If such a statement fails to be true, what  
are the implications for
the embedded curve $C$?}
\bigskip
\vskip10 truemm
\centerline{\bf BIBLIOGRAPHY}

\bigskip
\noindent
[A-C-G-H] Arbarello, E.- Cornalba, M. - Griffiths, P.A. - Harris, J.,
Geometry of Algebraic Curves I, Springer (1985).

\noindent
[A-S]\quad Arrondo, E. - Sols, I., On Congruences of Lines in the
Projective Space,  Soci\'et\'e Math\'ematique de France, M\'emoire n. 50
Suppl\'ement au Bulletin de la S.M.F., Tome 120, (1992), fascicule 3.

\noindent
[C-C-D]\quad  Chiantini, L. -  Ciliberto, C. -  Di Gennaro, V.,
The genus of projective curves, Duke Math. J., Vol. 70, No. 2, 
(1993), 229-245. 

\noindent
[C-K-M]\quad Clemens, H. - Koll\'ar, J. - Mori, S., Higher Dimensional
Complex Geometry, Ast\'erisque 166, Soci\'et\'e Math\'ematique 
de France (1988).

\noindent 
[E-P]\quad Ellinsgrud, G. - Peskine, C., Sur les surfaces lisses de
$\Bbb P^4$, Invent. Math., 95 (1989), 1-11.

\noindent
[G-P]\quad Gruson, L. - Peskine, C., Genre des courbes de 
l'espace projectif,
in Proceedings of Tromso, Conference on Algebraic Geometry. Springer
LNM 687 (1977), 31-59.

\noindent
[Ha]\quad Hartshorne, R., Algebraic Geometry, Springer GTM 52 (1977).

\noindent
[JH]\quad Harris, J., The Genus of Space Curves, Math. Ann. 249, 
191-204.
(1980)

\noindent
[Mi]\quad Migliore, J., Liaison of a union of Skew Lines in 
$\Bbb P^4$, Pac. Jour. Math. Vol 130, No.1, (1987) 153-170.

\noindent
[Mu]\quad Mumford, D., Lectures on Curves on an Algebraic Surface, 
Annals of
Mathematics Studies, n.59, Princeton Univ. Press (1966).

\enddocument